\documentclass{amsart}

\RequirePackage{fix-cm}
\usepackage{graphicx}
\usepackage{amsmath, amsopn,amstext,amscd,amsfonts,amssymb}
\usepackage{dsfont}
\usepackage{comment}
\usepackage[active]{srcltx}
\usepackage{graphicx, epsfig, subfig}

\usepackage{textcomp}

\usepackage{algorithm}

\makeatletter
\def\BState{\State\hskip-\ALG@thistlm}
\makeatother

\def\downbar#1{
\setbox10=\hbox{$#1$}
            \dimen10=\ht10 \advance\dimen10 by 2.5pt
            \ifdim \dimen10<15pt 
               \advance\dimen10 by -0.5pt
               \dimen11=\dimen10
               \advance\dimen10 by 2.5pt
               \lower \dimen11
            \else \lower \ht10 \fi
            \hbox {\hskip 1.5pt \vrule height \dimen10 depth \dp10}}
\def\upbar#1{
\setbox10=\hbox{$#1$}
            \dimen10=\ht10 \advance\dimen10 by \dp10 \advance\dimen10 by 2.5pt
            \ifdim \dimen10<15pt 
                \advance\dimen10 by 2pt \fi
            \raise 2.5pt \hbox {\hskip -1.5pt \vrule height \dimen10}}

\usepackage{multicol}
\usepackage{colortbl}
\usepackage{exscale}
\usepackage{amssymb,latexsym,amsthm,amsfonts,color,fancyhdr,mathrsfs}
\usepackage{lscape}
\usepackage{rotating}
\usepackage{caption}

\newtheorem{proposition}{\bf Proposition}[section]

\newtheorem{remark}{\bf Remark}[section]

\numberwithin{equation}{section}
\begin{document}

\title[On orthogonal polynomials described by Chebyshev polynomials]{A note on orthogonal polynomials described by Chebyshev polynomials}
\author{K. Castillo}
\address{CMUC, Department of Mathematics, University of Coimbra, 3001-501 Coimbra, Portugal}
\email{kenier@mat.uc.pt}
\author{M. N. de Jesus}
\address{CI$\&$DETS/IPV, Polytechnic Institute of Viseu, ESTGV, Campus Polit\'ecnico de Repeses, 3504-510 Viseu, Portugal}
\email{mnasce@estv.ipv.pt}
\author{J. Petronilho}
\address{CMUC, Department of Mathematics, University of Coimbra, 3001-501 Coimbra, Portugal}
\email{josep@mat.uc.pt}

\date{\today}

\begin{abstract}
The purpose of this note is to extend in a simple and unified way
some results on orthogonal polynomials with respect to the weight function
$$\frac{|T_m(x)|^p}{\sqrt{1-x^2}}\;,\quad-1<x<1\;,$$
where $T_m$ is the Chebyshev polynomial of the first kind of degree $m$ and $p>-1$.
\end{abstract}
\maketitle

\section{Main result}

Let $T_n$ and $U_n$ denote the Chebyshev polynomials of first and second kind,
that is, $T_n(\cos\theta)=\cos(n\theta)$ and $U_n(\cos\theta)=\sin\big((n+1)\theta\big)/\sin\theta$
for each nonnegative integer $n$ and $0<\theta<\pi$.
Let $\widehat{T}_n$ and $\widehat{U}_n$ denote the corresponding monic polynomials,
so that $\widehat{T}_0:=\widehat{U}_0:=1$, $\widehat{T}_n(x):=2^{1-n}T_n(x)$
and $\widehat{U}_n(x):=2^{-n}U_n(x)$ for each positive integer $n$.
Set $T_{-n}:=U_{-n}:=\widehat{T}_{-n}:=\widehat{U}_{-n}:=0$ for each positive integer $n$.
The reader is assumed familiar with basic properties of Chebyshev polynomials.
We prove the following

\begin{proposition}\label{main}
Fix an integer $m\geq2$. Define $t_0:=0$ and let $(t_n)_{n\geq1}$ be a sequence of nonzero complex numbers such that
\begin{align*}
& t_{2mn+j}=\mbox{$\frac14$}\;,\quad j\in\{0,1,\dots,2m-1\}\setminus\{0,1,m,m+1\}\;, \\
& t_{2mn}+t_{2mn+1}=t_{2mn+m}+t_{2mn+m+1}=\mbox{$\frac12$}\;,
\end{align*}
for each nonnegative integer $n$.
Let $(P_n)_{n\geq0}$ be the sequence of monic orthogonal polynomials given by
\begin{align}\label{pn}
P_{n+1}(x)=xP_n(x)-t_nP_{n-1}(x)\;.
\end{align}
Then $P_j(x)=\widehat{T}_{j}(x)$ for each $j\in\{0,1,\dots,m\}$ and
\begin{equation}\label{p2mn+m+j+1}
P_{2mn+m+j+1}(x) =\frac{A_j(n;x)Q_{n+1}\big(\widehat{T}_{2m}(x)\big)
+4^{-j}t_{2mn+m+1}B_j(n;x)Q_{n}\big(\widehat{T}_{2m}(x)\big)}{\widehat{U}_{m-1}(x)}
\end{equation}
for each $j\in\{0,1,\dots,2m-1\}$, where
\begin{align*}
A_j(n;x)&:=\widehat{U}_j(x)+\Big(\mbox{$\frac14$}-t_{2m(n+1)}\Big)
\Big(\widehat{U}_{j-m}(x)\widehat{U}_{m-2}(x)-\widehat{U}_{j-m-1}(x)\widehat{U}_{m-1}(x)\Big)\,, \\
B_j(n;x)&:=\widehat{U}_{2m-2-j}(x) \\
&\quad+\Big(\mbox{$\frac14$}-t_{2m(n+1)}\Big)
\Big(\widehat{U}_{m-j-3}(x)\widehat{U}_{m-1}(x)-\widehat{U}_{m-j-2}(x)\widehat{U}_{m-2}(x)\Big)\,,
\end{align*}
and $(Q_n)_{n\geq0}$ is the sequence of monic orthogonal polynomials given by 
\begin{align*}
Q_{n+1}(x)=(x-r_n)Q_n(x)-s_nQ_{n-1}(x)\;,
\end{align*}
with
\begin{align*}
r_n&:=\frac{1}{\displaystyle 2^{2m-4}}\Big(t_{2mn+m}t_{2mn+1}+t_{2m(n+1)}t_{2mn+m+1}-\mbox{$\frac18$}\Big)\;,\\
s_n&:=\frac{1}{\displaystyle 4^{2m-4}}\,t_{2mn}t_{2mn+1}t_{2mn+m}t_{2m(n-1)+m+1}\;.
\end{align*}
Assume furthermore that $t_n>0$ for each positive integer $n$.
Then $(P_n)_{n\geq0}$ and $(Q_n)_{n\geq0}$ are orthogonal polynomial sequences with respect to
certain positive measures, say $\mu_P$ and $\mu_Q$ respectively.
Suppose that $\mu_Q$ is absolutely continuous
with weight function $w_Q$ on $[\xi,\eta]$, with $-2^{1-2m}\leq\xi<\eta\leq2^{1-2m}$, i.e.,
$$
{\rm d}\mu_Q(x)=w_Q(x)\chi_{(\xi,\eta)}(x){\rm d}x\;.
$$
Suppose in addition that $w_Q$ satisfies the condition
$$C:=\int_\xi^\eta\frac{w_Q(x)}{x+2^{1-2m}}\,{\rm d}x<\infty\;.$$
Then $($up to a positive constant factor$)$
\begin{align}\label{muP}
{\rm d}\mu_P(x)=w_P(x)\chi_{E}(x){\rm d}x +
M\sum_{j=1}^m\delta\left(x-\cos\frac{(2j-1)\pi}{2m}\right)\;,
\end{align}
where
$$
w_P(x):=\left|\frac{U_{m-1}(x)}{T_m(x)}\right|w_Q\big(\widehat{T}_{2m}(x)\big)\;,\quad
x\in E:=\widehat{T}_{2m}^{-1}\big((\xi,\eta)\big)\,,
$$
and $M$ is a nonnegative number given by
$$
M:=\frac{1}{m}\left(\frac{2^{2m-3}\mu_Q(\mathbb{R})}{t_m}-C\right)\;.
$$
\end{proposition}

\begin{remark}
For $j=2m-1$, \eqref{p2mn+m+j+1} reduces to
$$
P_{2mn+m}(x)=\widehat{T}_m(x)\,Q_n\big(\widehat{T}_{2m}(x)\big)\;.
$$
\end{remark}

\begin{remark}\label{Zmzi}
In general, the set $E$ is the union of $2m$ disjoint open intervals
separated by the zeros of $U_{2m-1}$.
$($If all these intervals have the zeros of $U_{2m-1}$ as boundary points,
then $\overline{E}$ reduces to a single interval.$)$
Moreover,
$\mbox{\rm supp}\big(\mu_P\big)=\overline{E}$ if $M=0$, and
$\mbox{\rm supp}\big(\mu_P\big)=\overline{E}\cup Z_m$ if $M>0$,
where $Z_m$ is the set of zeros of $T_m$, i.e.,
$Z_m:=\big\{z_i:=\cos\big((2i-1)\pi/(2m)\big)\,|\, 1\leq i\leq m\big\}$.
\end{remark}

\begin{remark}
Proposition \ref{main} cannot be considered for $m=1$.
Indeed, if $m=1$ the conditions imposed on the sequence $(t_n)_{n\geq1}$ imply $t_2=0$.
In any case, after reading Section 2, interested readers are able to derive the corresponding
result for $m=1$.
\end{remark}

In order to illustrate the practical effectiveness of Proposition \ref{main},
we consider a simple example. 
Fix $p\in\mathbb{C}\setminus\{-1,-2,\ldots\}$ and an integer $m\geq2$. Define
$t_{2mn+j}:=1/4$ for each $j\in\{0,1,\ldots,2m-1\}\setminus\{0,1,m,m+1\}$ and
\begin{align*}
\qquad t_{2mn}&:=\frac{n}{4n+p}\;, &
t_{2mn+1}&:=\frac{2n+p}{2(4n+p)}\;, & \quad  \\
\qquad t_{2mn+m}&:=\frac{2n+p+1}{2(4n+p+2)}\;, &
t_{2mn+m+1}&:=\frac{2n+1}{2(4n+p+2)} &
\end{align*}
for each nonnegative integer $n$. (If $p=0$ it is understood that $t_1:=1/2$.)
These parameters satisfy the hypothesis of Proposition \ref{main}.
In this case 
\begin{align*}
r_n&=\frac{1}{2^{2m-1}}
\frac{p(p+2)}{(4n+p)(4n+p+4)}\;, \\
s_n&=\frac{1}{4^{2m-1}}\frac{2n(2n-1)(2n+p)(2n+p+1)}
{(4n+p-2)(4n+p)^2(4n+p+2)}\;,
\end{align*}
and so
$$Q_n(x)=2^{(1-2m)n}\widehat{P}_n^{\big(-\frac12,\frac{p+1}{2}\big)}\big(2^{2m-1}x\big)\;,$$
where $\widehat{P}_n^{(\alpha,\beta)}$ denotes the monic Jacobi polynomial of degree $n$
(see \cite[Chapter 4]{I2005} and \cite[(6.11)]{M1991}).
Hence, by Proposition \ref{main},
\begin{align}\label{ExamplePn}
&P_{2mn+m+j+1}(x) \\
&\quad=\frac{(2n+2+p)U_j(x)-pT_m(x)U_{j-m}(x)}{2^{(2m-1)n+m+j-1}(4n+p+4)U_{m-1}(x)}
\widehat{P}_{n+1}^{\big(-\frac12,\frac{p+1}{2}\big)}\big(T_{2m}(x)\big) \nonumber \\
&\quad\quad
+\frac{(2n+1)\big((2n+2)U_{2m-2-j}(x)-pU_{m-1}(x)T_{m-j-1}(x)\big)}
{2^{(2m-1)n+m+j}(4n+p+2)(4n+p+4)U_{m-1}(x)}
\widehat{P}_{n}^{\big(-\frac12,\frac{p+1}{2}\big)}\big(T_{2m}(x)\big) \nonumber
\end{align}
for each nonnegative integer $n$ and each $j\in\{0,1,\dots,2m-1\}$.
Furthermore, if $p>-1$ then $(Q_n)_{n\geq0}$ is a sequence of orthogonal polynomials associated with the weight function
$w_Q$ on $\big[-2^{1-2m},2^{1-2m}\big]$ given by
$$
w_Q(x)=2^{-p/2}\big(1-2^{2m-1}x\big)^{-\frac12}\big(1+2^{2m-1}x\big)^{\frac{p+1}{2}}\;,\quad -2^{1-2m}<x<2^{1-2m}\;.
$$
Since
\begin{align*}
\mu_Q(\mathbb{R})&=\int_{-2^{1-2m}}^{2^{1-2m}}w_Q(x)\,{\rm d}x
=2^{2-2m}\,\frac{p+1}{p+2}B\Big(\mbox{$\frac{p+1}{2}$},\mbox{$\frac{1}{2}$}\Big)\;,\\
C&=\int_{-2^{1-2m}}^{2^{1-2m}}\frac{w_Q(x)}{x+2^{1-2m}}\,{\rm d}x
=B\Big(\mbox{$\frac{p+1}{2}$},\mbox{$\frac{1}{2}$}\Big)\;,
\end{align*}
$B$ being the Beta function (see \cite[p. 8]{I2005}), we get $M=0$.
Taking into account that
$\overline{E}=\widehat{T}_{2m}^{-1}\big(\big[-2^{1-2m},2^{1-2m}\big]\big)={T}_{2m}^{-1}\big([-1,1]\big)=[-1,1]$,
we conclude that $(P_n)_{n\geq0}$ is a sequence of orthogonal polynomials associated with the weight function
$w_P$ on $[-1,1]$ given by
\begin{align}
w_P(x)&=2^{-p/2}\,\left|\frac{U_{m-1}(x)}{T_m(x)}\right|
\big(1-T_{2m}(x)\big)^{-\frac12}\big(1+T_{2m}(x)\big)^{\frac{p+1}{2}} \nonumber \\
&=\frac{\;|T_m(x)|^p}{\sqrt{1-x^2}}\;,\quad-1<x<1\;. \label{measureTmp}
\end{align}

\begin{remark}
The search of the recurrence coefficients from the weight function \eqref{measureTmp}
has arouse interest in recent and not recent years.
For $p=2$ see \cite{MVA1989};
for $p=1$ and $m=2$ see \cite{GLi1988};
for $p=2s$ with $s\in\mathbb{N}$ and $m=2$ see \cite{CMM2016};
and for $p=2s$ with $s>-1/2$ and $m=2$ see \cite{CMV2018}.
It is worth mentioning that in these works
the explicit representation \eqref{ExamplePn} has not been established.
While it is true that a representation given by a different polynomial mapping appears in \cite{MVA1989} for $p=2$
using results from \cite{GVA1988}.
\end{remark}

The sequence $(P_n)_{n\geq0}$ in Proposition \ref{main} may be regarded as an example of sieved orthogonal polynomials
$($see e.g. \cite{AAA1984,CI1986,CI1993,CIM1994}$)$. It is important to highlight that we can go even further and show more distinctly
how some recent developments of the theory of polynomial mappings (see \cite{MP2010,KMP2017,KMP2019}) apply to this
kind of problems, but that would demand a more extended discussion in which we
need not now get involved. Indeed, using the general results stated in \cite{KMP2017},
it can be shown that if $(P_n)_{n\geq0}$ or $(Q_n)_{n\geq0}$ in Proposition \ref{main} is semiclassical (see \cite{M1991}),
then so is the other one. In this case, we may obtain also the functional (distributional) equation
fulfilled by the (moment) regular functional associated with $(P_n)_{n\geq0}$ or $(Q_n)_{n\geq0}$.
In particular, as in \cite{KMP2019} or directly from \eqref{ExamplePn}, we may easily derive the linear homogeneous second order ordinary differential equation that the orthogonal polynomials with respect to \eqref{measureTmp} fulfil, which in turn leads to interesting electrostatic models (see \cite[Section 6]{KMP2019} and \cite[Section 3.5]{I2005}).

\section{Proof of Proposition \ref{main}}

Set $k:=2m$ and rewrite \eqref{pn} as a system of blocks of recurrence relations
\begin{align*}
x P_{nk+j}(x)=P_{nk+j+1}(x)+a_n^{(j)}P_{nk+j-1}(x) \;,\quad 0\leq j\leq k-1\;,
\end{align*}
where $P_{-1}:=0$ and $a_n^{(j)}:=t_{kn+j}$ whenever $(n,j)\neq(0,0)$.
Following \cite{CI1993,CIM1994}, we introduce the notation
\begin{equation}\nonumber
\Delta_n(i,j;x):=\left\{
\begin{array}{rl}
0\,,&j<i-2 \\
1\,,& j=i-2 \\
x\,,& j=i-1
\end{array}
\right.
\end{equation}
and
\begin{equation}\nonumber
\Delta_n(i,j;x):=\det\left(
\begin{array}{cccccc}
x & 1 & 0 &  \dots & 0 & 0  \\
a_n^{(i)} & x & 1 &  \dots & 0 & 0 \\
0 & a_n^{(i+1)} & x &   \dots & 0 & 0 \\
\vdots & \vdots & \vdots  & \ddots & \vdots & \vdots \\
0 & 0 & 0 &  \ldots & x & 1 \\
0 & 0 & 0 &  \ldots & a_n^{(j)} & x
\end{array}
\right)\;,\quad j\geq i\geq 1\;.
\end{equation}
The reader should satisfy himself that
$$
\Delta_n(m+2,m+j;x)=A_j(n;x)\,, \quad
\Delta_{n}(m+j+3,m+k-1;x)=B_j(n;x)
$$
for each $j\in\{0,1,\ldots,2m-1\}$. In particular,
\begin{align*}
\Delta_n(m+2,m+k-1;x)&=\widehat{U}_{2m-1}(x)=\widehat{T}_m(x)\widehat{U}_{m-1}(x) \\
&=\Delta_0(1,m-1;x)\eta_{k-1-m}(x)\,,\quad \eta_{k-1-m}(x):=\widehat{U}_{m-1}(x)\,.
\end{align*}
Moreover,
\begin{align*}
r_n(x)&:=a_{n}^{(m+1)}\Delta_{n}(m+3,m+k-1;x)-a_0^{(m+1)}\Delta_{0}(m+3,m+k-1;x) \\
&\quad+a_{n}^{(m)}\Delta_{n-1}(m+2,m+k-2;x)-a_0^{(m)}\Delta_{0}(1,m-2;x)\,\eta_{k-1-m}(x) \\
&=\frac{1}{\displaystyle 2^{2m-4}}\Big(t_{2mn+m}t_{2mn+1}+t_{2m(n+1)}t_{2mn+m+1}
-\mbox{$\frac12$}\,t_m-t_{m+1}t_{2m}\Big)
\end{align*}
for each positive integer $n$, and so
$$
\quad
r+r_n(0)=r_n\;,\quad r_0=r:=\frac{1}{\displaystyle 2^{2m-4}}\Big(\mbox{$\frac18$}-t_{m+1}t_{2m+1}\Big)\;. 
$$
Note also that
$a_{n}^{(m)}a_{n-1}^{(m+1)}\cdots a_{n-1}^{(m+k-1)}=s_n$.
In addition,
\begin{align*}
\pi_k(x)&:=\Delta_0(1,m;x)\,\eta_{k-1-m}(x)-a_0^{(m+1)}\,\Delta_0(m+3,m+k-1;x)+r\\
&=\widehat{T}_{2m}(x)\;.
\end{align*}
Consequently, the hypotheses of \cite[Theorem 2.1]{MP2010} are satisfied
with $\theta_m(x)=\widehat{T}_{m}(x)$, and so (\ref{p2mn+m+j+1}) follows.
Finally, the explicit representation of $\mu_P$ appearing in (\ref{muP}) follows by
\cite[Theorem 3.4\footnote{In \cite[Theorem 3.4]{MP2010}, $r=0$. Nevertheless, inspection
of its proof shows that the theorem remains true if instead $r=0$ we assume that the polynomials $\pi_k$ and $\theta_m\eta_{k-m-1}$ have real and distinct zeros.} and Remark 3.5]{MP2010}, after noting that
\begin{align*}
\pi_k'& =2m\widehat{U}_{2m-1}=2m\widehat{T}_{m}\widehat{U}_{m-1}=2m\theta_m\eta_{k-m-1}\;, \\
\pi_k(z_i)& =2^{1-2m}\,T_2\big(T_m(z_i)\big)=2^{1-2m}\,T_2(0)=-2^{1-2m}\;, \\
M_i & :=\frac{ \mu_Q(\mathbb{R})\,\Delta_0(2,m-1;z_i)/
\prod_{j=1}^m a_0^{(j)}-\eta_{k-1-m}(z_i)\,C}{ \theta_m^\prime(z_i)}\\
&=\frac{ \mu_Q(\mathbb{R})\,\widehat{U}_{m-1}(z_i)\cdot 2^{2m-3}t_m^{-1}
-\widehat{U}_{m-1}(z_i)\,C}{ m\widehat{U}_{m-1}(z_i)}=M
\end{align*}
for each $i\in\{1,\ldots,m\}$, where $z_i$ is a zero of $T_m$ (see Remark \ref{Zmzi}).

\section*{Acknowledgements}
KC and JP are supported by the Centre for Mathematics of the University of Coimbra --UID/MAT/00324/2019, funded by the Portuguese Government through FCT/MEC and co-funded by the European Regional Development Fund through the Partnership Agreement PT2020.
MNJ supported by UID/Multi/04016/2019, funded by FCT.
MNJ also thanks the Instituto Polit\'ecnico de Viseu and CI\&DETS for their support.

\end{document}